\documentclass{amsart}
\usepackage{amssymb,amscd,verbatim, amsthm,amsmath,amsgen,xspace}
\usepackage{latexsym}
\usepackage{amsfonts}
\usepackage{color}
\usepackage{ulem}
\usepackage{csquotes}
\usepackage[all]{xy}

\newcommand \unb[2]{\underset{#1}{{\underbrace{#2}}}}

\newcommand \dpg{\Delta^+}

\newcommand \dpk{\Delta^+_\frk}
\newcommand \dpl{\Delta^+_\frl}
\newcommand \du{\Delta(\fru)}

\newcommand \dg{\Delta}

\newcommand \wg{W_\frg}

\newcommand \wpg{W^1}
\newcommand \wl{W_\frl}
\newcommand \wlp{W_\frl^1}
\newcommand \wpl{W_\frl^1}
\newcommand \wk{W_\frk}

\newcommand \caH{\mathcal{H}}

\newcommand\ep{{\epsilon}}

\newcommand\la{{\lambda}}

\newcommand\Ad{\operatorname{Ad}}

\newcommand{\Kt}{\widetilde{K}}

\def\Ad{\mathop{\hbox {Ad}}\nolimits}

\DeclareMathOperator\im{Im}

\def\ker{\mathop{\hbox{Ker}}\nolimits}

\newcommand{\pf}{\begin{proof}}
\newcommand{\epf}{\end{proof}}
\newcommand{\eq}{\begin{equation}}
\newcommand{\eeq}{\end{equation}}
\newcommand{\eqn}{\begin{equation*}}
\newcommand{\eeqn}{\end{equation*}}

\newcommand{\aql}{A_\frq(\lambda)}
\newcommand{\aq}{A_\frq}

\newcommand{\frg}{\mathfrak{g}}

\newcommand{\frk}{\mathfrak{k}}
\newcommand{\frl}{\mathfrak{l}}

\newcommand{\fro}{\mathfrak{o}}
\newcommand{\frp}{\mathfrak{p}}
\newcommand{\frq}{\mathfrak{q}}

\newcommand{\frs}{\mathfrak{s}}
\newcommand{\frt}{\mathfrak{t}}
\newcommand{\fru}{\mathfrak{u}}

\newcommand{\frso}{\mathfrak{so}}
\newcommand{\frsp}{\mathfrak{sp}}

\newcommand{\bbC}{\mathbb{C}}

\newcommand{\bbN}{\mathbb{N}}

\newcommand{\bbR}{\mathbb{R}}

\newcommand{\bbZ}{\mathbb{Z}}
\newtheorem{thm}[equation]{Theorem}
\newtheorem{cor}[equation]{Corollary}
\newtheorem{lemma}[equation]{Lemma}
\newtheorem{prop}[equation]{Proposition}

\numberwithin{equation}{section}

\let\ssize\scriptstyle
\newif\ifFIRST\newdimen\MAXright\MAXright0pt
\def\sdynkin{\bgroup\eightpoint\dynkin}
\def\endsdynkin{\enddynkin\egroup}
\def\dynkin{\bgroup\FIRSTtrue\hskip.5em\setbox1\hbox{$\diagup$}%
	\setbox2\hbox{$\diagdown$}%
	\setbox0\hbox to2\wd1{\hrulefill}%
	\setbox3\hbox{$\bullet$}%
	\setbox4\hbox{$\times$}%
	\setbox7\hbox{$\circ$}
	\def\whiteroot##1{\ifFIRST\setbox5\hbox{$##1$}\ifdim\wd5>1.3em
		\hskip-.5em\hskip.5\wd5\fi\fi\FIRSTfalse
		\hskip-.25em\raise1.5\wd3\hbox to0pt{\hss\hskip.45em$
			\ssize##1$\hss}\copy7\hskip-.25em\setbox6\hbox{$##1$}
		\MAXright\wd6}
	\def\root##1{\ifFIRST\setbox5\hbox{$##1$}\ifdim\wd5>1.3em%
		\hskip-.5em\hskip.5\wd5\fi\fi\FIRSTfalse%
		\hskip-.25em\raise1.5\wd3\hbox to0pt{\hss\hskip.45em$%
			\ssize##1$\hss}\copy3\hskip-.25em\setbox6\hbox{$##1$}%
		\MAXright\wd6}%
	\def\whitedroot##1{\ifFIRST\setbox5\hbox{$##1$}\ifdim\wd5>1.3em
		\hskip-.5em\hskip.5\wd5\fi\fi\FIRSTfalse
		\hskip-.25em\lower1.8\wd3\hbox to0pt{\hss\hskip.45em$
			\ssize##1$\hss}\copy7\hskip-.25em\setbox6\hbox{$##1$}
		\MAXright\wd6}%
	\def\whiterroot##1{\hskip-.25em\copy7\hbox to0pt{\hskip.3em$\ssize##1$\hss}%
		\hskip-.25em\setbox6\hbox{\hskip.6em$##1##1$}%
		\MAXright\wd6}%
	\def\droot##1{\ifFIRST\setbox5\hbox{$##1$}\ifdim\wd5>1.3em%
		\hskip-.5em\hskip.5\wd5\fi\fi\FIRSTfalse%
		\hskip-.25em\lower1.8\wd3\hbox to0pt{\hss\hskip.45em$%
			\ssize##1$\hss}\copy3\hskip-.25em\setbox6\hbox{$##1$}%
		\MAXright\wd6}%
	\def\rroot##1{\hskip-.25em\copy3\hbox to0pt{\hskip.3em$\ssize##1$\hss}%
		\hskip-.25em\setbox6\hbox{\hskip.6em$##1##1$}%
		\MAXright\wd6}%
	\def\norroot##1{\hskip-.36em\copy4\hbox to0pt{\hskip.3em$\ssize##1$\hss}%
		\hskip-.48em\setbox6\hbox{\hskip.6em$##1##1$}%
		\MAXright\wd6}%
	\def\noroot##1{\ifFIRST\setbox5\hbox{$##1$}\ifdim\wd5>1.3em%
		\hskip-.5em\hskip.5\wd5\fi\fi\FIRSTfalse%
		\hskip-.36em\raise1.5\wd3\hbox to0pt{\hss\hskip.6em$%
			\ssize##1$\hss}\copy4\hskip-.38em\setbox6\hbox{$##1$}%
		\MAXright\wd6}%
	\def\nodroot##1{\ifFIRST\setbox5\hbox{$##1$}\ifdim\wd5>1.3em%
		\hskip-.5em\hskip.5\wd5\fi\fi\FIRSTfalse%
		\hskip-.36em\lower1.8\wd3\hbox to0pt{\hss\hskip.6em$%
			\ssize##1$\hss}\copy4\hskip-.38em\setbox6\hbox{$##1$}%
		\MAXright\wd6}%
	\def\nolink{\hskip\wd0}
	\def\link{\raise.22em\copy0}%
	\def\llink##1{\raise.32em\copy0\hskip-\wd0%
		\raise.12em\copy0\hskip-.5\wd0\hbox to0pt{\hss$##1$\hss}\hskip.5\wd0}%
	\def\lllink##1{\raise.22em\copy0\hskip-\wd0\raise.32em\copy0\hskip-\wd0%
		\raise.12em\copy0\hskip-.5\wd0\hbox to0pt{\hss$##1$\hss}\hskip.5\wd0}%
	\def\rootupright##1{\hbox to0pt{\raise.45em\copy1\hskip-.25em\raise1.3\ht1%
			\hbox{\copy3\hskip.3em$\ssize##1$}\hss}%
		\setbox6\hbox{\hskip.6em\copy1\copy1$##1##1$}%
		\ifdim\MAXright<\wd6\MAXright\wd6\fi}%
	\def\whiterootupright##1{\hbox to0pt{\raise.45em\copy1\hskip-.25em\raise1.3\ht1
			\hbox{\copy7\hskip.3em$\ssize##1$}\hss}
		\setbox6\hbox{\hskip.6em\copy1\copy1$##1##1$}
		\ifdim\MAXright<\wd6\MAXright\wd6\fi}
	\def\norootupright##1{\hbox to0pt{\raise.45em\copy1\hskip-.36em\raise1.3\ht1%
			\hbox{\copy4\hskip.3em$\ssize##1$}\hss}%
		\setbox6\hbox{\hskip.6em\copy1\copy1$##1##1$}%
		\ifdim\MAXright<\wd6\MAXright\wd6\fi}%
	\def\rootdownright##1{\hbox to0pt{\raise-.5em\copy2\hskip-.25em\raise-1.35\ht1%
			\hbox{\copy3\hskip.3em$\ssize##1$}\hss}\setbox6%
		\hbox{\hskip.6em\copy2\copy2$##1##1$}%
		\ifdim\MAXright<\wd6\MAXright\wd6\fi}%
	\def\whiterootdownright##1{\hbox to0pt{\raise-.5em\copy2\hskip-.25em\raise-1.35\ht1
			\hbox{\copy7\hskip.3em$\ssize##1$}\hss}\setbox6
		\hbox{\hskip.6em\copy2\copy2$##1##1$}
		\ifdim\MAXright<\wd6\MAXright\wd6\fi}
	\def\rootdown##1{\hbox to0pt{\hskip-.05em\vrule height.25em depth.65em%
			\hskip-.25em\raise-.95em\hbox{\copy3\hskip.3em$\ssize##1$}\hss}%
		\setbox6\hbox{$##1$}%
		\ifdim\MAXright<\wd6\MAXright\wd6\fi}%
	\def\whiterootdown##1{\hbox to0pt{\hskip-.05em\vrule height.25em depth.65em
			\hskip-.25em\raise-.95em\hbox{\copy7\hskip.3em$\ssize##1$}\hss}
		\setbox6\hbox{$##1$}
		\ifdim\MAXright<\wd6\MAXright\wd6\fi}
	\def\dots{\hskip.5em\cdots\hskip.5em}}%
\def\enddynkin{\ifdim\MAXright>1em\hskip.5\MAXright\else\hskip.5em\fi\egroup}%

\begin{document}

\bigskip
\title[Classification of $A_{q}(\lambda)$ modules]{Classification of $A_{q}(\lambda)$ modules by their Dirac cohomology for type $D$, $G_2$ and $\frsp(2n,\bbR)$}

\author{Ana Prli\'{c}}
\address{Department of Mathematics, University of Zagreb, Bijeni\v cka 30, 10000 Zagreb, Croatia}
\email{anaprlic@math.hr}
\thanks{The author is supported by grant no. 4176 of the Croatian Science Foundation and by the QuantiXLie Centre of Excellence, a project
cofinanced by the Croatian Government and European Union through the
European Regional Development Fund - the Competitiveness and Cohesion
Operational Programme (Grant KK.01.1.1.01.0004). }
	
\begin{abstract} Let $G$ be a connected real reductive group with maximal compact subgroup $K$ of the same rank as $G$. In the recent paper of Huang, Pand\v{z}i\'{c} and Vogan, it was shown that the admissible $\Theta$--stable parabolic subalgebras $\frq$ of $\frg$ are in one-to-one correspodence with the faces of $W \rho$ intersecting the $\frk$--dominant Weyl chamber and that $A_{\frq}(0)$--modules can be classified by their Dirac cohomology in geometric terms. They described in detail the cases when $\frg_0$ is of type $A$, $B$, $F$ and $C$ except for $\frg_0 = \frsp(2n, \bbR)$.  We will describe faces corresponding to $A_{\frq}(0)$--modules for $\frg_0 = \frsp(2n, \bbR)$ and for $\frg_0$ of type $D$ and $G_2$.
\end{abstract}

\maketitle

\section{Introduction}\label{section intro}
Let $G$ be a connected real reductive group with Cartan involution $\Theta$, such that $K = G^{\Theta}$ is a maximal compact subgroup of $G$, and let 
\[
\frg_0 = \frk_0 + \frp_0
\]
be the Cartan decomposition of the Lie algebra of $G$ corresponding to $\Theta$. We will assume that $\frg_0$ and $\frk_0$ have equal rank. Let $\frt_0$ be the common Cartan subalgebra of $\frg_0$ and $\frk_0$.

The $A_{\frq}(\lambda)$--modules were introduced by Vogan and Zuckerman in \cite{VZ}. It was shown in \cite{SR1} that any irreducible unitary $(\frg, K)$--modules with strongly regular infinitesimal character is isomorphic to some $A_{\frq}(\lambda)$. The $A_{q}(\lambda)$--modules are also important for applications to automorphic forms. 

The notion of the Dirac cohomology was introduced by Vogan in \cite{V2}. Let $U(\frg)$ denote the universal enveloping algebra of $\frg$, and let $C(\frp)$ denote the Clifford algebra of $\frp$ with respect to the Killing form $B$. Let $b_i$ be any basis of $\frp$, and let $d_i$ be the dual basis of $\frp$ with respect to $B$. The Dirac operator $D\in U(\frg)\otimes C(\frp)$ is defined as
\[
D=\sum_i b_i\otimes d_i.
\]

The Dirac operator is independent of the choice of $b_i$, and $K$-invariant with respect to $\Ad\otimes\Ad$.

Let $X$ be a $(\frg,K)$-module and let $S$ be a spin module for $C(\frp)$. Then $X\otimes S$ is $(U(\frg)\otimes C(\frp), \tilde{K})$--module, where $\tilde{K}$ is the spin double cover of $K$. The Dirac cohomology of $X$ is
\[
H_D(X)=\ker D/(\im D\cap\ker D).
\]
It is a $\Kt$-module, which is finite-dimensional in case $X$ is of finite length. If $X$ is unitary, then there is an inner product on $X\otimes S$ such that $D$ is self adjoint. It follows that
\[
H_D(X)=\ker D=\ker D^2.
\]

The main result about Dirac cohomology is the following theorem conjectured by Vogan \cite{V2} and proved by Huang and Pand\v{z}i\'{c} (\cite{HP1}). For a dominant $\frk$-weight $\gamma\in\frt^*$, we denote by $E_\gamma$ the irreducible finite-dimensional $\frk$-module with infinitesimal character $\gamma$. 

\begin{thm}{\bf (\cite{HP1})}
Let $X$ be an irreducible $(\frg,K)$-module with nonzero Dirac cohomology. Let $E_\gamma$
be a $\widetilde{K}$-type contained in $H_D(X)$. Then the infinitesimal character of $X$ is
conjugate to $\gamma$ under the Weyl group $W$ of $\frg$ with respect to $\frt$.
\end{thm}
For more details about Dirac cohomology and its properties see \cite{HP2}. In \cite{HPV}, $A_{\frq}(0)$--modules (shortly, $A_{\frq}$--modules) are classifed by their Dirac cohomology in very explicit, geometric terms. It was shown that if the Levi subalgebra $\frl_0$ has no compact simple factors and no factors isomorphic to $\frso(2n, 1)$, where $n > 1$, or $\frsp(p, q)$, where $p, q > 0$, or the nonsplit $F_4$, then the $A_{\frq}$--module is uniquely determined by its Dirac cohomology. This is the case for $\frg_0$ of type $A, D, E$ or $G$. Most of the remaining simple concompact equal rank real Lie algebras (i.e. Lie algebras of type $B, C$ and $F$), do contain Levi subalgebas with factors isomorphic to $\frso(2n, 1)$, where $n > 1$, or $\frsp(p, q)$, where $p, q > 0$ and $A_{\frq}$--modules are not uniquely determined by their Dirac cohomology. However, in \cite{HPV} it is listed which modules share the same Dirac cohomology. 

Recall that in \cite{KV}, the $\Theta$--stable parabolic subalgebras of $\frg_0$ are constructed via $\frk$--dominant elements $H$ of $i \frt_0$. We will also use this construction. Let us just mention that there are several ways of constructing $\frq$ and the connection between some of them is described in \cite[Section~4.]{HPV}.

As in \cite{HPV}, we concentrate only on the case $\lambda = 0$. By \cite[Proposition~5.2.]{HPV}, the admissible $\frq$ are in one-to-one correspodence with the faces of $W \rho$ intersecting the $\frk$--dominant Weyl chamber $C_{\frk}$, and the corresponding face $\Phi$ is equal $W_{\frl} \rho'$, where $\rho'$ is any element of the face $\Phi$. The Dirac cohomology of $A_{\frq}$ corresponds to $\Phi \cap C_{\frk}$.

We will describe faces corresponding to $\Theta$--stable parabolic subalgebras for the cases $\frg_0 = \frsp(2n, \bbR)$ and for $\frg_0$ of type $D$ and $G_2$. The rest of the cases for classical equal rank $\frg_0$, as well as the case when $\frg_0$ is of type $F$, are described in \cite{HPV}. A description of faces corresponding to $\Theta$--stable parabolic subalgebras for the case when $\frg_0$ is of type $E$ turned out to be much complicated and too long and we omit it. That case will be considered separately in future research and probably requires the aid of a computer.

The paper is organised as follows. In Section 2, we give a brief review of the main results about $A_{\frq}(\lambda)$--modules, Dirac cohomology and their description in geometric terms, as it was done in \cite{HPV}. In Section 3 we describe the faces corresponding to admissible $\Theta$--stable parabolic subalgebras of $\frsp(2n, \bbR)$, in Section 4 we describe the faces corresponding to admissible $\Theta$--stable parabolic subalgebras of $\frso(2p, 2q)$, in Section 5 we describe the faces corresponding to admissible $\Theta$--stable parabolic subalgebras of $\frso^{*}(2n)$ and in Section 6 we describe the faces corresponding to admissible $\Theta$--stable parabolic subalgebras of $\frg_0$ of type $G_2$.

In future, we hope to generalize our results to the unequal rank Lie algebras. The author would like to thank Pavle Pand\v{z}i\'{c} for valuable comments which helped to improve the paper.

\section{$A_{q}(\lambda)$--modules, Dirac cohomology and faces of $W \rho$}\label{section prelim}
Let us recall the definition of $\Theta$--stable parabolic subalgebras from \cite{KV}. With the notation as in the introduction, let $H \in i \frt_0$. Since the elements from $\text{ad}\frt_0$ are skew symmetric, all roots have real values on $H$. The parabolic subalgebra $\frq$ of $\frg$ corresponding to $H$ is 
\[
\frq = \frl \oplus \fru, 
\]
where $\frl$ is the kernel of $\text{ad}_{\frg}H$, and $\fru$ is the sum of the eigenspaces of $\text{ad}_{\frg}H$ for positive eigenvalues. In terms of the root system $\Delta(\frg, \frt)$,
\[
\frl = \frt \oplus \bigoplus_{\alpha \in \Delta_{\frg}, \alpha(H) = 0} \frg_{\alpha}; \qquad \fru = \bigoplus_{\alpha \in \Delta_{\frg}, \alpha(H) > 0} \frg_{\alpha},
\]
where $\frg_\alpha$ denotes the root space corresponding to the root $\alpha$. In particular, $\Delta_\frl$ can be identified with the subset of $\dg$ consisting of roots that are 0 on $H$.

We will always assume that the element $H$ defining $\frq$ is $\frk$-dominant, i.e.,
\[
\alpha(H)\geq 0,\quad \text{for every}\quad \alpha\in\Delta^+_\frk.
\]
Here $\Delta_{\frk}^{+}$ is a positive root system for $\frk$ which we fix throughout the paper. 	
 
As in \cite{SR1}, we call the parabolic subalgebras corresponding to such $H$ admissible.

The Levi subgroup $L$ of $G$ corresponding to $\frq$ is defined as
\[
L=N_G(\frq)=\{g\in G\,\big|\, \Ad(g)\frq=\frq\}.
\]

As in \cite{VZ}, $\lambda\in\frt^*$ is admissible for $\frq$ if $\lambda$ is the differential of a unitary character $\bbC_\lambda$ of $L$, and satisfies the condition
\[
\langle\lambda,\alpha\rangle\geq 0,\qquad \alpha\in\Delta(\fru).
\]
We assume that $\lambda$ is $\frk$-dominant. Let us recall a definition of $A_\frq(\lambda)$ modules.

\begin{thm}{\bf (\cite{VZ},\cite{V1})}
Let $\frq=\frl\oplus\fru$ be a $\theta$-stable parabolic subalgebra of $\frg$ and let $\lambda\in\frt^*$ be admissible as
defined above. Fix a choice of $\dpl$ containing $\Delta_{\frl \cap \frk}^{+}$ and set $\dpg=\dpl\cup\du$.
Then there is a unique irreducible unitary $(\frg,K)$-module $A_\frq(\lambda)$ with the following properties:

(i) $A_\frq(\lambda)$ has infinitesimal character $\lambda+\rho$;

(ii) The lowest $K$-type of $A_\frq(\lambda)$ has highest weight
$\mu(\frq,\lambda)=\lambda+2\rho(\fru\cap\frp)$;

(iii) All $K$-types of $A_\frq(\lambda)$ have highest weights of the form
\eqn
\mu(\frq,\lambda)+\sum_{\beta\in\Delta(\fru\cap\frp)}n_\beta\beta
\eeqn
with $n_\beta$ non-negative integers.
\end{thm}

The $A_\frq(\lambda)$ modules are constructed by cohomological induction described in \cite{KV} and in \cite{HP2}. They are irreducible and unitary. They are important since they \textquote{cover} all irreducible unitary $(\frg,K)$-modules with strongly regular infinitesimal character (\cite{SR1}), and also because they show up in many applications. The next theorem tells us when two $A_{\frq}(\lambda)$--modules are equal.

\begin{prop}{\bf (\cite[Proposition~4.5.]{SR2})}
\label{propSR}
Let $\frq$ and $\frq'$ be two admissible $\theta$-stable parabolic subalgebras of $\frg$ as above, and let $\lambda$ respectively $\lambda'$ be admissible for $\frq$ respectively $\frq'$. Then
\[
A_{\frq}(\lambda)\cong A_{\frq'}(\lambda')
\]
if and only if $\lambda=\lambda'$ and
\[
\fru\cap\frp=\fru'\cap\frp.
\]
\end{prop}

We will consider only $\frq$ for which $\frl_0$ has no compact simple factors (see \cite[Corollary~2.4.]{HPV}).

Now we recall the description of the Dirac cohomology of $A_\frq(\lambda)$. Let us fix a positive root system $\dpg\supset\dpk$. We denote by $C_\frg\subset C_\frk$ the corresponding dominant Weyl chambers for $\frg$ respectively $\frk$.

\begin{lemma}{\bf (\cite[Lemma~3.2.]{HPV})}
\label{wpg}
For an element $w$ of $W$, the following statements are equivalent:
\begin{enumerate}
\item $wC_\frg\subset C_\frk$;
\item $w\rho\in C_\frk$, where $\rho$ denotes the half sum of roots in $\dpg$;
\item The positive root system $w\dpg$ contains $\Delta_\frk^+$;
\item $w$ is the unique shortest element in its right $\wk$-coset $\wk w$, where $W_{\frk}$ denotes the Weyl group of $\frk$ with respect to $\frt$.
\end{enumerate}
\end{lemma}

The set of all $w\in W$ satisfying the equivalent conditions of Lemma \ref{wpg} is denoted by $\wpg$. 

\begin{thm}[\cite{HKP}]
\label{HKP}
Let $\frq$ be an admissible $\theta$-stable parabolic subalgebra of $\frg$, and let $\lambda\in\frt^*$ be admissible for $\frq$.
Then the Dirac cohomology of the unitary irreducible $(\frg,K)$-module
$A_\frq(\lambda)$ is
\[
H_D(A_\frq(\lambda))=\bigoplus_{w\in W_\frl^1}E_{w(\lambda+\rho)},
\]
where $W_\frl^1$ is the subset of the Weyl group $W_\frl$ of $\frl$, defined analogously to $\wpg$. 
\end{thm}
It follows from \cite[Lemma~4.1.]{HPV} that $W_\frl^1=\wpg\cap W_\frl$. For an $\aql$ module as above, we can choose positive root systems
$\Delta_\frl^+$ and $\dpg=\Delta_\frl^+\cup\Delta(\fru)$, such that $\la+\rho$ is dominant with respect to $\dpg$.

As in \cite{HPV}, we can concentrate on the case $\la=0$, i.e., the case of modules with infinitesimal character equal to $\rho$. In the rest of the paper we study only the modules $A_\frq=A_\frq(0)$. 

If we identify $H_D(\aq)$ with the set of $\frk$-infinitesimal characters appearing in it, it becomes exactly the set of all $\rho$ for the positive root systems $\dpg\supset\dpk$ compatible with $\frq$ (i.e., containing $\Delta(\fru)$). Fixing one such $\dpg=\dpl\cup\du$ and the corresponding $\rho$, and defining $\wpl$ with respect to $\dpl$, this set becomes
\eq
\label{face}
\wpl\rho=\{w\rho\,\big|\, w\in \wlp\}.
\eeq
Let $P_\rho=C(W\rho)$ denote the polytope obtained as the convex hull of $W\rho$. A a face of $P_\rho$ is either all of $P_\rho$, or a nonempty intersection of $P_\rho$ with a hyperplane $\caH$ in $\frt^*_\bbR$, such that all of $P_\rho$ is on one side of $\caH$. We will denote by $\Phi\subseteq W\rho$ the set of vertices of the face, and the face itself will be denoted by $C(\Phi)$, the convex hull of $\Phi$.

For a fixed face $C(\Phi)$, we consider a vector $N$ at a point of $C(\Phi)$, orthogonal to $C(\Phi)$ and pointing away from $P_{\rho}$, and define the corresponding $H \in i \frt_0$ by
\eq
\label{NH}
\langle N,\mu\rangle = \mu(H),\qquad \mu\in\frt^*_\bbR.
\eeq
Thus $\mu\in W\rho$ is in $\Phi$ if and only if $\mu(H)\geq\nu(H)$ for all $\nu\in W\rho$. We will call both $\Phi$ and $C(\Phi)$ the face corresponding to $H$, and we will use the notation $\Phi=\Phi_H$.

\begin{lemma}{\bf (\cite[Lemma~5.1.]{HPV})}
\label{kdom}
The face $\Phi_H$ intersects the $\frk$-dominant chamber $C_\frk$ precisely when the corresponding $H\in i\frt_0$ is $\frk$-dominant, i.e.,
\[
\alpha(H) \geq 0,\qquad \alpha\in\dpk.
\]
\end{lemma}

\begin{thm}{\bf (\cite[Theorem~5.2.]{HPV})}
\label{q and faces}
Admissible $\theta$-stable parabolic subalgebras $\frq=\frl\oplus\fru$ of $\frg$ are in one-to-one correspondence with the faces of $W\rho$ intersecting the $\frk$-dominant chamber $C_\frk$. If $\Phi=\Phi_\frq$ corresponds to $\frq$ under this correspondence, then $\Phi=\wl\rho'$, where $\rho'$ is any element of $\Phi$.
\end{thm}

If the face $\Phi=\Phi_\frq$ corresponds to $\frq$, then $\Phi\cap C_\frk$ is the Dirac cohomology of the module $\aq$.

The algorithm for getting the face $\Phi_H$ from a given $\frk$--dominant $H \in i\frt_0$ is as follows:
First we find the shortest $w \in W^{1}$ such that $w^{-1}H$ is $\frg$--dominant. Then we define $\rho' = w \rho$, where $\rho$ is the half-sum of positive roots of the positive fixed root system for $\frg$ (containing our fixed $\Delta_{\frk}^{+}$). Then the face $\Phi_H$ equals $W_{\frl} \rho'$. The face $\Phi_H$ corresponds to $A_{\frq}$ (where $\frq$ is defined via $H$), while $\Phi_H \cap C_{\frk}$ corresponds to the Dirac cohomology of $A_{\frq}$.

\section{Classification of $A_{q}$ modules for $sp(2n,\bbR)$}\label{sec:sp2n}

Let $\frg_0=\frsp(2n,\bbR)$. In this case,
$\frk_0=\fru(n)$. In the standard positive root system, the compact positive roots are 
\[
\epsilon_i-\epsilon_j,\qquad 1\leq i<j\leq n,
\]
while the noncompact positive roots are 
\[
\epsilon_i+\epsilon_j,\quad 1\leq i<j\leq n\quad\text{ and }\ 2\ep_i,\quad 1\leq i\leq n. 
\]
The standard simple roots are $\epsilon_1-\epsilon_2,\dots,\epsilon_{n-1}-\epsilon_n,2\epsilon_n$, with the only noncompact one being $2\ep_n$. The standard diagram thus has one black point at the last place $n$.

We will denote the elements of $i\frt_0$ or of $\frt^*_\bbR$, which are both isomorphic to $\bbC^n$, by $(x_1,x_2\dots,x_n)$.
Such an element is $\frk$-dominant if and only if
\[
x_1\geq \dots\geq x_n.
\]
The set $\wpg$ consists of $w\in W$ which take
\[
\rho=(n,n-1,\dots,2,1)
\]
to $(i_1,\dots,i_n)$, with $i_1>i_2>\dots>i_n$ such that $|i_1|,\dots,|i_n|$ are all the coordinates of $\rho$. 

The admissible $\theta$-stable parabolic subalgebras with no compact factors are described by the following proposition.

\begin{prop}
\label{levis_sp2}
Let $\frg_0=\frs\frp(2n,\bbR)$ with $n\geq 2$. Then 
\begin{enumerate}
\item[(i)]
The maximal (proper) Levi subalgebras of $\frg_0$ with no compact factors are 
\[
\frl_0^{p,q}=\fru(p,q)\times \frsp(2r,\bbR),
\]
where $p,q\in\bbZ^+$ are such that $p+q\leq n$, 
\[
p=0\,\Rightarrow\,q=1;\qquad  q=0\,\Rightarrow\, p=1
\]
and $r= n - p - q$.
\item[(ii)]
Each $\frl_0^{p,q}$ as in (i) is embedded into $\frg_0$ in a unique way: the factor $\fru(p,q)$ has simple roots
\[
\ep_1-\ep_2,\dots,\ep_{p-1}-\ep_p,\ep_p+\ep_{p+r+1},-\ep_{p+r+1}+\ep_{p+r+2},\dots,-\ep_{n-1}+\ep_n,
\]
while the factor $\frsp(2r,\bbR)$ is built on coordinates $p+1,\dots,p+r$.

If $p\neq q$, then $\frl_0^{p,q}\cong\frl_0^{q,p}$, but the embeddings of these two Levi subalgebras are different. There are no other isomorphisms between the $\frl_0^{p,q}$.

\item[(iii)] 
Each of the maximal Levi subalgebras $\frl^{p,q}$ described above is the Levi subalgebra of a unique admissible $\theta$-stable parabolic subalgebra $\frq^{p,q}=\frl^{p,q}\oplus \fru^{p,q}$ of $\frg$.
The corresponding $H\in i\frt_0$ are
\[
H^{p,q}=(\unb{p}{1,\dots,1},\unb{r}{0,\dots,0},\unb{q}{-1,\dots,-1}).
\]
The roots of $\fru^{p,q}$ can easily be written down: they are the roots that are positive on $H^{p,q}$.

\item[(iv)] Let $\frq=\frl\oplus\fru$ be any admissible $\theta$-stable parabolic subalgebra of $\frg$, such that $\frl_0$ has no compact factors. Then 
\[
\frl_0=\frl'_0\times\frsp(2r,\bbR)\subseteq \frl_0^{p,q},
\] 
with $\frl_0^{p,q}$ one of the maximal Levi subalgebras described above, and with $\frl_0'$ a Levi subalgebra of the first factor $\fru(p,q)$ of $\frl_0^{p,q}$. All possible $\frl'_0$ are described in \cite[Proposition~6.1.]{HPV}.

Moreover, there is a parabolic subalgebra $\frq_1=\frl_0'\oplus\fru_1$ of $\fru(p,q)$, such that 
\[
\fru=\fru_1\oplus\fru^{p,q}.
\]
\end{enumerate}
\end{prop}
\pf 
It is clear that the $H^{p,q}$ given in (iii) define admissible $\theta$-stable parabolic subalgebras with Levi subalgebras as in the statement of the proposition. We have to show that there are no other possibilities for a maximal Levi subalgebra.

Let $\frl$ be any Levi subalgebra of $\frg$ without compact factors defined by some $\frk$--dominant $H$. Without loss of generality we may assume that $H$ has integer coordinates. Since $H$ is $\frk$--dominant, $H$ is of the form
\[
H = (\unb{p_l}{l,\dots,l}, \unb{p_{l-1}}{l-1,\dots,l-1}, \dots \unb{p_1}{1,\dots,1}, \unb{r}{0,\dots,0}, \unb{q_1}{-1,\dots,-1}, \dots, \unb{q_l}{-l,\dots,-l}),
\]
where $p_i, q_j \in \mathbb{Z}^{+}$. Since $\frl$ has no compact simple factors, then 
\[
p_i=0\ \Rightarrow\ q_i \leq 1;\qquad q_j=0\ \Rightarrow\ p_j \leq 1.
\]
Of course, it is possible that for some $i$ both $p_i = q_i = 0$, but the corresponding Levi subalgebra can also be constructed via some $H$ as the above where $(p_k, q_k) \neq (0,0)$ for all $k$, so we can assume 
\[
p_i=0\ \Rightarrow\ q_i = 1;\qquad q_j=0\ \Rightarrow\ p_j = 1.
\]
Obviously, 
\[
H^{'} = (\unb{p}{1,\dots,1}, \unb{r}{0,\dots,0}, \unb{q}{-1,\dots,-1}),
\] 
where $p = p_1 + p_2 + \cdots + p_l$, $q = q_1 + q_2 + \cdots + q_l$ defines a strictly larger $\frl$.

This proves (i) and (ii).

(iii) and (iv) are also proved, since we have seen that there are no further possibilities for $H$ corresponding to Levi subalgebras without compact factors besides the ones listed; in particular, each $\frl$ corresponds to only one $H$ and hence to only one $\frq$. 
\epf

Putting the results in Proposition \ref{levis_sp2} together with \cite[Proposition~6.1.]{HPV} , we obtain

\begin{cor}
\label{combparC2}
Let $\frg_0=\frsp(2n,\bbR)$ with $n\geq 2$. 
The modules $A_\frq$, or equivalently the admissible $\theta$-stable parabolic subalgebras $\frq=\frl\oplus\fru$ of $\frg$ such that $\frl_0$ has no compact factors, are in one-to-one correspondence with ordered sequences
\[
(p_1,q_1),(p_2,q_2),\dots,(p_l,q_l),
\]
with $l,p_i,q_j\in\bbZ^+$, such that
\begin{gather*}
p_1+p_2+\dots+p_l+q_1+q_2+\dots+q_l\leq n;\\
\text{and } \qquad p_i=0\ \Rightarrow\ q_i=1;\qquad q_j=0\ \Rightarrow\ p_j=1.
\end{gather*}
\qed
\end{cor}

Let us describe the faces of $W \rho$ corresponding to the $A_{\frq}$ modules given by Proposition \ref{levis_sp2}. 

The shortest $w\in\wpg$ such that $w^{-1}H^{p,q}$ is $\frg$-dominant with respect to the standard $\dpg$ sends 
$\rho=(n,n-1,\dots,1)$ to
\[
\rho'=(n,\dots,n-p+1;r,\dots,1;-(r+1),\dots,-(n-p)).
\]
The vertices of the corresponding face $\Phi = W_{\frl}\rho'$ are obtained from $\rho'$ by:

\begin{enumerate}
\item permuting the first $p$ coordinates, permuting the last $q$ coordinates, exchanging pairs $(a, -b)$ and $(b, -a)$ with the first component among the first
$p$ coordinates and the second component among the last $q$ coordinates;
\item permuting and making arbitrary sign changes to the group of coordinates
in between.
\end{enumerate}

Similarly, the shortest $w\in\wpg$ such that 
\[
w^{-1}H = w^{-1}(\unb{p_l}{l,\dots,l}, \unb{p_{l-1}}{l-1,\dots,l-1}, \dots \unb{p_1}{1,\dots,1}, \unb{r}{0,\dots,0}, \unb{q_1}{-1,\dots,-1}, \dots, \unb{q_l}{-l,\dots,-l})
\]
is $\frg$-dominant with respect to the standard $\dpg$ sends
$\rho=(n,n-1,\dots,1)$ to
\begin{align*}
\rho''=(& \unb{p_l}{n, n-1, \dots,n - p_l + 1}, \unb{p_{l-1}}{n - (p_l + q_l),\dots,n - (p_l + q_l) - p_{l-1} + 1}, \dots, \\
& \unb{p_1}{r + p_1 + q_1,\dots,r + q_1 + 1}, \unb{r}{r, r-1,\dots,1}, \unb{q_1}{-(r + 1), - (r +
2), \dots,-(r + q_1)}, \\
& \unb{q_2}{ - (r + p_1 + q_1 + 1), \dots,  - (r + p_1 + q_1 + q_2)}, \dots, \unb{q_l}{ - n + p_l + q_l - 1,\dots. - n + p_l}).
\end{align*}

The vertices of the corresponding face $\Phi = W_{\frl} \rho''$ are obtained from $\rho''$ by:
\begin{enumerate}
\item permuting coordinates in groups that correspond to repeated coordinates of $H$.
\item exchanging pairs $(a, -b)$ and $(b, -a)$ in groups of coordinates with the first component in the group that corresponds to repeated coordinates $t, \cdots, t$ of $H$, and the second component in the group that corresponds to repeated coordinates $-t, \cdots, -t$ of $H$, where $t \in \{1, \cdots, l\}$.
\item permuting and making arbitrary sign changes to the group of coordinates at the places $p+1, \cdots p+r$.
\end{enumerate}

Now $\Phi\cap C_\frk$, which corresponds to $H_D(\aq)$ for the $\frq$ corresponding to $\Phi$, consists of the vertices of $\Phi$ that are $\frk$-dominant, i.e., the vertices with coordinates in descending order. The descending order can be checked separately for each group of coordinates. From this one can see that each $\Phi$ is determined by $\Phi\cap C_\frk$, so the modules $\aq$ are determined by their Dirac cohomology. This also follows from \cite[Section~5.]{HPV}, since Levi subalgebras of $\frg_0$ can not contain factors $\frso(2n,1)$, $\frsp(p,q)$, or the nonsplit $F_4$. The last statement can be easily seen by choosing $\frl_0$ from a Vogan diagram. If the long root is chosen, it is always going to be black (noncompact), and the three bad cases are exactly those with all long roots compact.

We summarize the results for the case $\frg_0=\frsp(2n,\bbR)$ in the following theorem.

\begin{thm}
\label{thmsp2n}
If $\frg_0=\frsp(2n,\bbR)$, then the modules $\aq$ correspond to faces $\Phi$ of $W\rho$ intersecting $C_\frk$, which can be described by the following procedure:

\begin{enumerate}
\item Pick a $\frk$-dominant conjugate $\rho'$ of $\rho$;
\item Pick an ordered sequence $(p_1,q_1),\dots,(p_l,q_l)$ as in Corollary \ref{combparC2}; let $p=p_1+\dots+p_l$ and $q=q_1+\dots+q_l$;
\item Consider the following groups of coordinates of $\rho'$:
\begin{multline*}
\rho'_1,\dots,\rho'_{p_l};\rho'_{n-q_l+1},\dots,\rho'_{n};\qquad \rho'_{p_l+1},\dots,\rho'_{p_l+p_{l-1}};\rho'_{n-q_l - q_{l-1} + 1},\dots,\rho'_{n-q_l};\qquad\dots \\
\dots\qquad \rho'_{p-p_1+1},\dots,\rho'_{p};\rho'_{n - q +1},\dots,\rho'_{n - q + q_1}.
\end{multline*}
Each of these groups of coordinates has two subgroups, separated by a semicolon. In the first of these subgroups, all coordinates are positive, and in the second subgroup the coordinates are negative. Now we first ignore the signs, permute the whole group in all possible ways, and then put back the signs to the coordinates in the new second subgroup.

The coordinates that are not part of the above groups of coordinates form the last group 
\[
\rho'_{p + 1},\dots,\rho'_{p+r}.
\]

To these coordinates we apply all permutations and arbitrary sign changes.
\end{enumerate}
\item If an $\aq$ module corresponds to a face $\Phi$, then
its Dirac cohomology $H_D(A_\frq)$ consists of the $\Kt$-types with infinitesimal characters in $\Phi\cap C_\frk$.
$H_D(\aq)$, or $\Phi\cap C_\frk$, uniquely determines the module $A_\frq$.  The set $\Phi\cap C_\frk$ consists of vertices that have decreasing coordinates, and this property can be checked for the whole vertex at once, or group by group. 
\qed
\end{thm}

\section{Classification of $A_{q}$ modules for $\frso(2p, 2q)$}\label{sec:D1}

We consider $\frg_0=\frs\fro(2p,2q)$, $p,q>0$, $p+q \geq 4$. The cases when $p+q<4$ are covered by other types and we omit them. 

For $\frg_0=\frso(2p,2q)$, we have 
$\frk_0=\frso(2p)\times\frso(2q)$. In the standard positive root system, the compact positive roots are 
\[
\epsilon_i\pm\epsilon_j,\quad 1\leq i<j\leq p\ \text{ or }\ p+1\leq i<j\leq p+q.
\]
The noncompact positive roots are 
\[
\epsilon_i\pm\epsilon_j,\quad 1\leq i\leq p\ \text{ and }\ p+1\leq j\leq p+q. 
\]
We will denote the elements of $i\frt_0$ or of $\frt^*_\bbR$, which are both isomorphic to $\bbC^{p+q}\cong\bbC^p\times\bbC^q$, by
\[
(x_1,x_2\dots,x_p\,|\, y_1,y_2,\dots,y_q).
\]
Such an element is $\frk$-dominant if and only if
\[
x_1\geq x_2\geq \dots\geq |x_p|\quad\text{and}\quad y_1\geq y_2\geq\dots\geq |y_q|.
\]

The set $W_\frg^1$ corresponding to the standard positive system consists of $w\in W_\frg$ which take
\[
\rho=(n-1,n-2,\dots,1,0),
\]
where $n = p + q$, to $(i_1,\dots,i_p,j_1,\dots,j_q)$, with $i_1>i_2>\dots>i_{p-1}>|i_p|$ and $j_1>j_2>\dots>|j_q|$ being (all) coordinates of $\rho$. Thus, $w$ is a $(p,q)$-shuffle, followed by a possible sign change at places $p$ and $p+q$. Note that one of these coordinates is 0, so the sign change will not have effect. However, elements of $\wg$ must have an even number of sign changes, so we need to have both sign changes if any.

We now turn to describing the admissible $\theta$-stable parabolic subalgebras of $\frg$.

\begin{prop} 
\label{leviD} Let $\frg_0=\frso(2p,2q)$, with $p,q>0$ and $p+q \geq 4$. Then any maximal proper Levi subalgebra $\frl_0$ of $\frg_0$ without compact factors is one of the following:
\begin{enumerate}
\item $\frl_0=\fru(p-r,q-s)\times\frso(2r,2s)$ for some $r$ and $s$ with $0< r\leq p$ and $0< s\leq q$, such that $r=p$ implies $s=q-1$, and $s=q$ implies $r=p-1$. This Levi subalgebra is embedded in a unique way, with the $\fru(p-r,q-s)$ factor having roots
\[
\ep_i-\ep_j,\qquad 1\leq i\leq p-r;\quad p+1\leq j\leq p+q-s,
\]
and with the $\frso(2r,2s)$ being built on coordinates
$p-r+1,\dots,p$ and $p+q-s+1,\dots,p+q$.

The corresponding $H\in i\frt_0$ are
\[
H=(\unb{p-r}{1,\dots,1},\unb{r}{0,\dots,0}\,|\,\unb{q-s}{1,\dots,1},\unb{s}{0,\dots,0}),
\]
The roots of $\fru^{p,q}$ can easily be written down: they are the roots that are positive on $H$.

\item $\frl_0=\fru(p,q)$, embedded in one of the following four ways:
\begin{enumerate}
\item[(a)] $\frl_0$ has simple roots
\[
\ep_1-\ep_2,\dots,\ep_{n-1}-\ep_n,
\]
where $n=p+q$. The corresponding $H\in i\frt_0$ is
\[
H_1=(1,\dots,1\,|\, 1,\dots,1).
\]
\item[(b)] $\frl_0$ has simple roots
\[
\ep_1-\ep_2,\dots,\ep_{p-2}-\ep_{p-1},\ep_{p-1}+\ep_p, -\ep_p-\ep_{p+1},\ep_{p+1}-\ep_{p+2},\dots,\ep_{n-1}-\ep_n.
\]
The corresponding $H\in i\frt_0$ is
\[
H_2=(1,\dots,1, -1\,|\, 1,\dots,1).
\]
\item[(c)] $\frl_0$ has simple roots
\[
\ep_1-\ep_2,\dots,\ep_{n-2}-\ep_{n-1},\ep_{n-1}+\ep_n.
\]
The corresponding $H\in i\frt_0$ is
\[
H_3=(1,\dots,1\,|\, 1,\dots,1, -1).
\]
\item[(d)] $\frl_0$ has simple roots
\[
\ep_1-\ep_2,\dots,\ep_{p-2}-\ep_{p-1},\ep_{p-1}+\ep_p, -\ep_p-\ep_{p+1},\ep_{p+1}-\ep_{p+2},\dots,\ep_{n-2}-\ep_{n-1},\ep_{n-1}+\ep_n.
\]
The corresponding $H\in i\frt_0$ is
\[
H_4=(1,\dots,1, -1\,|\, 1,\dots,1, -1).
\]
\end{enumerate}
\end{enumerate}

Any Levi subalgebra of $\frg_0$ with no compact factors is either of the form $\frl_0'\times\frso(2r,2s)$ with $\frl_0'$ a Levi subalgebra of $\fru(p-r,q-s)$ described in (1), or it is a Levi subalgebra of one of the four copies of $\fru(p,q)$ described in (2).

\end{prop}
\pf 
Let $\frl_0$ be any Levi subalgebra of $\frso(2p, 2q)$ without compact factors. As in the case of $\frg = \frsp(2n, \bbR)$, we can assume that all coordinates of the corresponding $H$ are integers. One posibillity is that all coordinates of $H$ are nonnegative. Since $H$ is $\frk$--dominant, $H$ is of the form
\begin{align*}
H_{r, s} =( & \unb{p_l}{l,\dots,l}, \unb{p_{l-1}}{l-1,\dots,l-1}, \dots, \unb{p_1}{1,\dots,1}, \unb{r}{0,\dots,0} \,|\, \\
&  \unb{q_l}{l,\dots,l}, \unb{q_{l-1}}{l-1,\dots,l-1}, \dots \unb{q_1}{1,\dots,1}, \unb{s}{0,\dots,0}),
\end{align*}
where $l \in \bbN$.
It is clear that $\frl_0$ is contained in the Levi subalgebra corresponding to 
\[
H=(\unb{p-r}{1,\dots,1},\unb{r}{0,\dots,0}\,|\,\unb{q-s}{1,\dots,1},\unb{s}{0,\dots,0}).
\]
If $r = 0$, then  $s \leq 1$ and if $s = 0$, $r \leq 1$, because, otherwise, $\frl_0$ would have a compact factor. If $(r, s) \in \{(0, 0), (0, 1), (1, 0) \}$, the Levi subalgebra $\frl_0$ is contained in the Levi subalgebra corresponding to 
\[
H_1=(1,\dots,1 \,|\, 1,\dots,1 ).
\]
If $H$ contains negative coordinates, then, since $H$ is $\frk$--dominant, those coordinates can be only at the places $p$ or $p + q$ and $H$ can not contain more than one zero coordinate. Otherwise either $H$ is not $\frk$--dominant or the corresponding $\frq$ contains a simple factor. This means that $H$ is of the form
\[
H=(\unb{p_l}{l,\dots,l}, \unb{p_{l-1}}{l-1,\dots,l-1}, \dots, \unb{p_1}{1,\dots,1}, p_0 \,|\, \unb{q_l}{l,\dots,l}, \unb{q_{l-1}}{l-1,\dots,l-1}, \dots, \unb{q_1}{1,\dots,1}, q_0),
\]
where at least one of the coordinates $p_0$ and $q_0$ is negative. If $p_0 < 0$ and $q_0 \geq 0$, then $\frl_0$ is contained in the Levi subalgebra corresponding to  
\[
H_2=(1,\dots,1, -1 \,|\, 1,\dots,1 ).
\]
If $p_0 \geq 0$ and $q_0 < 0$, then $\frl_0$ is contained in the Levi subalgebra corresponding to  
\[
H_3=(1,\dots,1 \,|\, 1,\dots,1, -1 ).
\]
Finally, if both $p_0$ and $q_0$ are negative, then $\frl_0$ is contained in the Levi subalgebra corresponding to  
\[
H_4=(1,\dots,1, -1 \,|\, 1,\dots,1, -1 ).
\]
This finishes the proof. 
\epf

Putting the results in Proposition \ref{leviD} together with \cite[Proposition~6.1.]{HPV}, we obtain
\begin{cor}
\label{combparD}
Let $\frg_0=\frso(2p,2q)$ where  $p, q > 0$ and $p + q \geq 4$. 
The modules $A_\frq$, or equivalently the admissible $\theta$-stable parabolic subalgebras $\frq=\frl\oplus\fru$ of $\frg$ such that $\frl_0$ has no compact factors, are in one-to-one correspondence with ordered sequences
\[
(p_1,q_1),(p_2,q_2),\dots,(p_l,q_l), (r, s), (n_p, n_q) 
\]
with $l,p_i,q_j, r, s \in\bbZ^+$, $n_p, n_q \in \{ 0, 1 \}$ such that
\begin{gather*}
p_1+p_2+\dots+p_l + r = p;\\
q_1+q_2+\dots+q_l + s = q; \\
p_i=0\ \Rightarrow\ q_i=1;\qquad q_j=0\ \Rightarrow\ p_j=1; \\
r=p \Rightarrow\ s=q-1; \qquad  s=q \Rightarrow\ r=p-1; \\
r = 0 \Rightarrow\ s \leq 1; \qquad s= 0 \Rightarrow\ r \leq 1; \\
r > 0 \Rightarrow\ n_p = 0; \qquad s > 0 \Rightarrow\ n_q = 0;\\
n_p = 1 \Rightarrow\  r = 0,  p_1 \geq 1; \qquad n_q = 1 \Rightarrow\  s = 0, q_1 \geq 1.
\end{gather*}
\qed
\end{cor}

\pf
If $r > 0$ and $s > 0$, then $(n_p, n_q) = (0, 0)$. From Proposition \ref{leviD} and \cite[Proposition~6.1.]{HPV} it follows that the parabolic subalgebra with parameters
\[
(p_1,q_1),(p_2,q_2),\dots,(p_l,q_l), (r, s), (0, 0) 
\]
has Levi subalgebra $\frl_0$ isomorphic to
\[
\fru(p_1, q_1) \times \dots \times \fru(p_l, q_l) \times \frso(2r, 2s). 
\]
If $(r, s) = (0, 0)$ and $(n_p, n_q) = (0, 0)$, then the parabolic subalgebra with parameters
\[
(p_1,q_1),(p_2,q_2),\dots,(p_l,q_l), (0, 0), (0, 0) 
\]
has Levi subalgebra $\frl_0$ isomorphic to
\[
\fru(p_1, q_1) \times \dots \times \fru(p_l, q_l). 
\]
with the corresponding $H$ equal to 
\[
H^{0, 0} =(\unb{p_l}{l,\dots,l}, \unb{p_{l-1}}{l-1,\dots,l-1}, \dots, \unb{p_1}{1,\dots,1} \,|\, \unb{q_l}{l,\dots,l}, \unb{q_{l-1}}{l-1,\dots,l-1}, \dots, \unb{q_1}{1,\dots,1}).
\]
If $(n_p, n_q) = (1, 0)$, then $r = 0$, $s \leq 1$ and $p_1 > 0$. The parabolic subalgebra with parameters
\[
(p_1,q_1),(p_2,q_2),\dots,(p_l,q_l), (0, s), (1, 0) 
\]
has Levi subalgebra $\frl_0$ isomorphic to
\[
\fru(p_1, q_1) \times \dots \times \fru(p_l, q_l),
\]
with the corresponding $H$ equal to 
\begin{align*}
H^{1, 0, s} =(& \unb{p_l}{l,\dots,l}, \unb{p_{l-1}}{l-1,\dots,l-1}, \dots, \unb{p_1}{1,\dots,1, -1} \,|\, \\
& \unb{q_l}{l,\dots,l}, \unb{q_{l-1}}{l-1,\dots,l-1}, \dots, \unb{q_1}{1,\dots,1}, \unb{s}{0}).
\end{align*}
If $(n_p, n_q) = (0, 1)$, then $s = 0$, $r \leq 1$ and $q_1 > 0$. The parabolic subalgebra with parameters
\[
(p_1,q_1),(p_2,q_2),\dots,(p_l,q_l), (r, 0), (0, 1) 
\]
has Levi subalgebra $\frl_0$ isomorphic to
\[
\fru(p_1, q_1) \times \dots \times \fru(p_l, q_l),
\]
with the corresponding $H$ equal to 
\begin{align*}
H^{0, 1, r} =(& \unb{p_l}{l,\dots,l}, \unb{p_{l-1}}{l-1,\dots,l-1}, \dots, \unb{p_1}{1,\dots,1} \unb{r}{0}\,|\, \\
& \unb{q_l}{l,\dots,l}, \unb{q_{l-1}}{l-1,\dots,l-1}, \dots, \unb{q_1}{1,\dots,1, -1}).
\end{align*}
Finally, if $(n_p, n_q) = (1, 1)$, then $(r,s) = (0, 0)$, $p_1 > 0$, $q_1 > 0$. The parabolic subalgebra with parameters, 
\[
(p_1,q_1),(p_2,q_2),\dots,(p_l,q_l), (0, 0), (1, 1) 
\]
has Levi subalgebra $\frl_0$ isomorphic to
\[
\fru(p_1, q_1) \times \dots \times \fru(p_l, q_l),
\]
with the corresponding $H$ equal to 
\begin{align*}
H^{1, 1} =(& \unb{p_l}{l,\dots,l}, \unb{p_{l-1}}{l-1,\dots,l-1}, \dots, \unb{p_1}{1,\dots,1, -1} \,|\, \\
& \unb{q_l}{l,\dots,l}, \unb{q_{l-1}}{l-1,\dots,l-1}, \dots, \unb{q_1}{1,\dots,1, -1}).
\end{align*}
\epf

Let us dscribe the faces of $W \rho$ corresponding to $A_{\frq}$ modules given by the proposition \ref{leviD}. 

The shortest $w\in\wpg$ such that
\begin{align*}
w^{-1} H_{r,s} = w^{-1}(& \unb{p_l}{l,\dots,l}, \unb{p_{l-1}}{l-1,\dots,l-1}, \dots, \unb{p_1}{1,\dots,1}, \unb{r}{0,\dots,0} \,|\, \\
& \unb{q_l}{l,\dots,l}, \unb{q_{l-1}}{l-1,\dots,l-1}, \dots, \unb{q_1}{1,\dots,1}, \unb{s}{0,\dots,0})
\end{align*}
is $\frg$-dominant with respect to the standard $\dpg$ sends 
$\rho=(n - 1, n - 2,\dots,1, 0)$ to
\begin{align*}
\rho'  =(& \unb{p_l}{n-1, n-2, \dots,n - p_l}, \unb{p_{l-1}}{n - (p_l + q_l) - 1,\dots,n - (p_l + q_l) - p_{l-1}}, \dots \\
& \unb{p_1}{r + s + p_1 + q_1 - 1,\dots,r + s + q_1}, \unb{r}{r + s - 1, \dots, s}, \unb{q_l}{n - p_l - 1, \dots, n - (p_l + q_l)}, \\
& \unb{q_{l-1}}{ n - (p_l + q_l) - p_{l-1} - 1, \dots,  n - (p_l + q_l) - (p_{l-1} + q_{l-1})}, \dots, \unb{q_1}{r + s - 1 + q_1, \dots, r + s}, \\
& \unb{s}{s - 1, \dots, 1, 0}).
\end{align*}

The vertices of the corresponding face $\Phi = W_{\frl} \rho'$ are obtained from $\rho'$ by:
\begin{enumerate}
\item permuting coordinates in groups that correspond to repeated coordinates of $H_{r, s}$.

\item exchanging pairs $(a, b)$ and $(-b, -a)$ in the group of coordinates at the places $p, p - 1, \cdots, p - r + 1$ and $p + q, p + q - 1, \cdots p + q - s + 1$.
\end{enumerate}

The shortest $w\in\wpg$ such that 
\begin{align*}
w^{-1}H^{0, 0} =w^{-1}(& \unb{p_l}{l,\dots,l}, \unb{p_{l-1}}{l-1,\dots,l-1}, \dots, \unb{p_1}{1,\dots,1} \,|\, \\
& \unb{q_l}{l,\dots,l}, \unb{q_{l-1}}{l-1,\dots,l-1}, \dots, \unb{q_1}{1,\dots,1}).
\end{align*}
is $\frg$-dominant with respect to the standard $\dpg$ sends 
$\rho=(n - 1, n - 2,\dots,1, 0)$ to
\begin{align*}
\rho'  =(& \unb{p_l}{n-1, n-2, \dots,n - p_l}, \unb{p_{l-1}}{n - (p_l + q_l) - 1,\dots,n - (p_l + q_l) - p_{l-1}}, \dots, \\
& \unb{p_1}{p_1 + q_1 - 1,\dots, q_1}, \unb{q_l}{n - p_l - 1, \dots, n - (p_l + q_l)}, \\
& \unb{q_{l-1}}{ n - (p_l + q_l) - p_{l-1} - 1, \dots,  n - (p_l + q_l) - (p_{l-1} + q_{l-1})}, \dots, \unb{q_1}{q_1 - 1, \dots, 0}).
\end{align*}

The vertices of the corresponding face $\Phi = W_{\frl} \rho'$ are obtained from $\rho'$ by permuting coordinates in groups that correspond to repeated coordinates of $H$.

The shortest $w\in\wpg$ such that 
\begin{align*}
w^{-1}H^{1, 0, s} =w^{-1}(& \unb{p_l}{l,\dots,l}, \unb{p_{l-1}}{l-1,\dots,l-1}, \dots, \unb{p_1}{1,\dots,1, -1} \,|\, \\
& \unb{q_l}{l,\dots,l}, \unb{q_{l-1}}{l-1,\dots,l-1}, \dots, \unb{q_1}{1,\dots,1}, \unb{s}{0})
\end{align*}
is $\frg$-dominant with respect to the standard $\dpg$ sends 
$\rho=(n - 1, n - 2,\dots,1, 0)$ to
\begin{align*}
\rho'  =(& \unb{p_l}{n-1, n-2, \dots,n - p_l}, \unb{p_{l-1}}{n - (p_l + q_l) - 1,\dots,n - (p_l + q_l) - p_{l-1}}, \dots, \\
& \unb{p_1 - 1}{p_1 + q_1 - 1, \dots, q_1 + 1}, 0, \unb{q_l}{n - p_l - 1, \dots, n - (p_l + q_l)}, \\
& \unb{q_{l-1}}{ n - (p_l + q_l) - p_{l-1} - 1, \dots,  n - (p_l + q_l) - (p_{l-1} + q_{l-1})}, \dots, \unb{q_1}{q_1, \dots, 2, 1})
\end{align*}
if $s = 0$ and to 
\begin{align*}
\rho'  =(& \unb{p_l}{n-1, n-2, \dots,n - p_l}, \unb{p_{l-1}}{n - (p_l + q_l) - 1,\dots,n - (p_l + q_l) - p_{l-1}}, \dots, \\
& \unb{p_1 - 1}{p_1 + q_1 - 1, \dots, q_1 + 1}, -1, \unb{q_l}{n - p_l - 1, \dots, n - (p_l + q_l)}, \\
& \unb{q_{l-1}}{ n - (p_l + q_l) - p_{l-1} - 1, \dots,  n - (p_l + q_l) - (p_{l-1} + q_{l-1})}, \dots, \unb{q_1}{q_1, \dots, 2, 0})
\end{align*}
if $s = 1$.

The vertices of the corresponding face $\Phi = W_{\frl} \rho'$ are obtained from $\rho'$ by:
\begin{enumerate}
\item permuting coordinates in groups that correspond to repeated coordinates of $H$.

\item exchanging pairs $(a, b)$ and $(-b, -a)$, where one of $a$ or $b$ is the $p$th coordinate of $\rho'$, and the second one is in the group of coordinates of $\rho'$ at the places $p - 1, \cdots, p - p_1 + 1$ and $p + q - s, p + q - 1 - s, \cdots, p + q - q_1 + 1 - s$.
\end{enumerate}

The shortest $w\in\wpg$ such that 
\begin{align*}
w^{-1}H^{1, 0, r} =w^{-1}(& \unb{p_l}{l,\dots,l}, \unb{p_{l-1}}{l-1,\dots,l-1}, \dots, \unb{p_1}{1,\dots,1}, \unb{r}{0}\,|\, \\
& \unb{q_l}{l,\dots,l}, \unb{q_{l-1}}{l-1,\dots,l-1}, \dots, \unb{q_1}{1,\dots,1, -1})
\end{align*}
is $\frg$-dominant with respect to the standard $\dpg$ sends 
$\rho=(n - 1, n - 2,\dots,1, 0)$ to
\begin{align*}
\rho'  =(& \unb{p_l}{n-1, n-2, \dots,n - p_l}, \unb{p_{l-1}}{n - (p_l + q_l) - 1,\dots,n - (p_l + q_l) - p_{l-1}}, \dots, \\
& \unb{p_1}{p_1 + q_1 - 1,\dots, q_1}, \unb{q_l}{n - p_l - 1, \dots, n - (p_l + q_l)}, \\
& \unb{q_{l-1}}{ n - (p_l + q_l) - p_{l-1} - 1, \dots,  n - (p_l + q_l) - (p_{l-1} + q_{l-1})}, \dots, \unb{q_1}{q_1 - 1, \dots, 1, 0})
\end{align*}
if  $r = 0$, and to 
\begin{align*}
\rho'  =(& \unb{p_l}{n-1, n-2, \dots,n - p_l}, \unb{p_{l-1}}{n - (p_l + q_l) - 1,\dots,n - (p_l + q_l) - p_{l-1}}, \dots, \\
& \unb{p_1}{p_1 + q_1,\dots, q_1 + 1}, 0, \unb{q_l}{n - p_l - 1, \dots, n - (p_l + q_l)}, \\
& \unb{q_{l-1}}{ n - (p_l + q_l) - p_{l-1} - 1, \dots,  n - (p_l + q_l) - (p_{l-1} + q_{l-1})}, \dots, \unb{q_1}{q_1, \dots, 2, -1})
\end{align*}
if $r = 1$.

The vertices of the corresponding face $\Phi = W_{\frl} \rho'$ are obtained from $\rho'$ by:
\begin{enumerate}
\item permuting coordinates in groups that correspond to repeated coordinates of $H$.

\item exchanging pairs $(a, b)$ and $(-b, -a)$, where one of $a$ or $b$ is the $p + q$th coordinate of $\rho'$, and the second one is in the group of coordinates of $\rho'$ at the places $p - r, p - r - 1 \cdots, p - p_1 + 1 - r$ and $p + q - 1, p + q - 2, \cdots p + q - q_1 + 1$.
\end{enumerate}

The shortest $w\in\wpg$ such that 
\begin{align*}
w^{-1}H^{1, 1} =w^{-1}(& \unb{p_l}{l,\dots,l}, \unb{p_{l-1}}{l-1,\dots,l-1}, \dots, \unb{p_1 - 1}{1,\dots,1}, -1 \,|\, \\
& \unb{q_l}{l,\dots,l}, \unb{q_{l-1}}{l-1,\dots,l-1}, \dots, \unb{q_1}{1,\dots,1}, -1)
\end{align*}
is $\frg$-dominant with respect to the standard $\dpg$ sends 
$\rho=(n - 1, n - 2,\dots,1, 0)$ to
\begin{align*}
\rho' =(& \unb{p_l}{n-1, n-2, \dots,n - p_l}, \unb{p_{l-1}}{n - (p_l + q_l) - 1,\dots,n - (p_l + q_l) - p_{l-1}}, \dots, \\
& \unb{p_1 - 1}{p_1 + q_1 - 1,\dots, q_1 + 1}, 0, \unb{q_l}{n - p_l - 1, \dots, n - (p_l + q_l)}, \\
& \unb{q_{l-1}}{ n - (p_l + q_l) - p_{l-1} - 1, \dots,  n - (p_l + q_l) - (p_{l-1} + q_{l-1})}, \dots, \unb{q_1 - 1}{q_1, \dots, 2}, -1).
\end{align*}

The vertices of the corresponding face $\Phi = W_{\frl} \rho'$ are obtained from $\rho'$ by:
\begin{enumerate}
\item permuting coordinates in groups that correspond to repeated coordinates of $H$.

\item exchanging pairs $(a, b)$ and $(-b, -a)$, where one of $a$ or $b$ is the $p$th or $p + q$th  coordinate of $\rho'$, and the second one is in the group of coordinates of $\rho'$ at the places $p - 1, p - 2, \cdots, p - p_1 + 1$ and $p + q - 1, p + q - 2, \cdots p + q - q_1 + 1$.
\end{enumerate}

We summarize the results for the case $\frg_0=\frso(2p, 2q)$, $p, q > 0, p + q \geq 4$ in the following theorem.

\begin{thm}
\label{thmso2}
If $\frg_0=\frso(2p,2q)$, then the modules $\aq$ correspond to faces $\Phi$ of $W\rho$ intersecting $C_\frk$, which can be described by the following procedure:

\begin{enumerate}
\item Pick a $\frk$-dominant conjugate $\rho'$ of $\rho$;
\item Pick an ordered sequence $(p_1,q_1),\dots,(p_l,q_l), (r,s), (n_p, n_q)$ as in Corollary \ref{combparD}
\item 
\begin{enumerate}
\item If $(n_p, n_q) = (0, 0)$, then the vertices of $\Phi$ are obtained by permutations of the following groups of coordinates of $\rho'$:
\begin{multline*}
\rho'_1,\dots,\rho'_{p_l};\rho'_{p + 1},\dots,\rho'_{p + q_l};\qquad \rho'_{p_l + 1},\dots,\rho'_{p_l + p_{l-1}};\rho'_{p + q_l + 1},\dots,\rho'_{p + q_l + q_{l-1}};\qquad \dots \\
\dots\qquad \rho'_{p - r - p_1 + 1},\dots,\rho'_{p - r};\rho'_{p + q - s - q_1 + 1},\dots,\rho'_{p + q - s}.
\end{multline*}

The coordinates that are not part of the above groups of coordinates form the last group 
\[
\rho'_{p - r + 1},\dots,\rho'_{p}; \rho'_{p + q - s + 1},\dots,\rho'_{p + q}
\]

To these coordinates we apply all permutations and even number of sign changes.

\item If $(n_p, n_q) = (1, 0)$, then the $p$th coordinate of $\rho'$ is negative or zero. We get the vertices of $\Phi$ by ignoring the signs, permuting the following groups of coordinates of $\rho'$:
\begin{multline*}
\rho'_1,\dots,\rho'_{p_l};\rho'_{p + 1},\dots,\rho'_{p + q_l};\qquad \rho'_{p_l + 1},\dots,\rho'_{p_l + p_{l-1}};\rho'_{p + q_l + 1},\dots,\rho'_{p + q_l + q_{l-1}};\qquad \dots \\
\dots\qquad \rho'_{p - p_1 + 1},\dots,\rho'_{p};\rho'_{p + q - q_1 + 1 - s},\dots,\rho'_{p + q - s},
\end{multline*}
and then putting the negative sign to the $p$th coordinate.  

\item If $(n_p, n_q) = (0, 1)$, then the $(p + q)$th coordinate of $\rho'$ is negative or zero. We get the vertices of $\Phi$ by ignoring the signs, permuting the following groups of coordinates of $\rho'$:
\begin{multline*}
\rho'_1,\dots,\rho'_{p_l};\rho'_{p + 1},\dots,\rho'_{p + q_l};\qquad \rho'_{p_l + 1},\dots,\rho'_{p_l + p_{l-1}};\rho'_{p + q_l + 1},\dots,\rho'_{p + q_l + q_{l-1}};\qquad \dots \\
\dots\qquad \rho'_{p - p_1 + 1 - r},\dots,\rho'_{p - r};\rho'_{p + q - q_1 + 1},\dots,\rho'_{p + q},
\end{multline*}
and then putting the negative sign to the $(p + q)$th coordinate.  

\item If $(n_p, n_q) = (1, 1)$, then $p$th and $(p + q)$th coordinates of $\rho'$ are negative or zero. We get the vertices of $\Phi$ by ignoring the signs, permuting the following groups of coordinates of $\rho'$:
\begin{multline*}
\rho'_1,\dots,\rho'_{p_l};\rho'_{p + 1},\dots,\rho'_{p + q_l};\qquad \rho'_{p_l + 1},\dots,\rho'_{p_l + p_{l-1}};\rho'_{p + q_l + 1},\dots,\rho'_{p + q_l + q_{l-1}};\qquad \dots \\
\dots\qquad \rho'_{p - p_1 + 1},\dots,\rho'_{p};\rho'_{p + q - q_1 + 1},\dots,\rho'_{p + q},
\end{multline*}
and then putting the negative sign to the $p$th and $(p + q)$th coordinate.  
\end{enumerate}
\end{enumerate}
\item If an $\aq$ module corresponds to a face $\Phi$, then
its Dirac cohomology $H_D(A_\frq)$ consists of the $\Kt$-types with infinitesimal characters in $\Phi\cap C_\frk$.
$H_D(\aq)$, or $\Phi\cap C_\frk$, uniquely determines the module $A_\frq$.  
\qed
\end{thm}

\section{Classification of $A_{q}$ modules for $\frso^{*}(2n)$}\label{sec:D2}

Let us now consider the case $\frg_0=\frso^*(2n)$, $n \geq 4$.
In this case, $\frk_0=\fru(n)$. In the standard positive root system, the compact positive roots are 
\[
\epsilon_i-\epsilon_j,\quad 1\leq i<j\leq n,
\]
and the noncompact positive roots are 
\[
\epsilon_i+\epsilon_j,\quad 1\leq i<j\leq n,
\]
We will denote the elements of $i\frt_0$ or of $\frt^*_\bbR$, which are both isomorphic to $\bbC^n$, by $(x_1,x_2\dots,x_n)$.
Such an element is $\frk$-dominant if and only if
\[
x_1\geq x_2\geq \dots\geq x_n.
\]
The set $W_\frg^1$ corresponding to the standard positive system consists of $w\in W_\frg$ which take
\[
\rho=(n - 1, n - 2, \dots, 1, 0)
\]
to $(i_1,\dots,i_n)$, with $i_1>i_2>\dots>i_n$ and and such that $|i_1|,\dots,|i_n|$ are all coordinates of $\rho$. Such $w$ can be given by a choice of signs $\eta_1,\dots,\eta_n$, with -1 appearing an even number of times. After putting these signs on the respective coordinates of $\rho$, there is a unique permutation bringing the result to the $\frk$-dominant chamber.

We now turn to describing the admissible $\theta$-stable parabolic subalgebras of $\frg$. The classification is similar to the case $\frsp(2n,\bbR)$.

\begin{prop}
\label{leviD2}
Let $\frg_0=\frso^*(2n)$ with $n\geq4$. Then the maximal (proper) Levi subalgebras of $\frg_0$ with no compact factors are given as follows. Let $p,q$ and $r$ be nonnegative integers satisfying 
\[
p+q+r=n;\qquad p=0\,\Rightarrow\,q=1;\qquad  q=0\,\Rightarrow\, p=1.
\]
For each such choice of $p,q,r$ there is a maximal Levi subalgebra
\[
\frl_0=\fru(p,q)\times \frso^*(2r)
\]
of $\frg_0$ embedded in the following way:
The factor $\fru(p,q)$ has simple roots
\[
\ep_1-\ep_2,\dots,\ep_{p-1}-\ep_p,\ep_p+\ep_{p+r+1},-\ep_{p+r+1}+\ep_{p+r+2},\dots,-\ep_{n-1}+\ep_n,
\]
while the factor $\frso^*(2r)$ is built on coordinates $p+1,\dots,p+r$.

If $p\neq q$, then the Levi subalgebra corresponding to $p,q,r$ is isomorphic to the Levi subalgebra corresponding to $q,p,r$, but the embeddings of these two Levi subalgebras are different.

Any Levi subalgebra $\frl_0$ of $\frg_0$ without compact factors is of the form 
\[
\frl_0=\frl^{'}_0\times\frso^*(2r),
\] 
such that $\frl_0^{'}$ is a Levi subalgebra of $\fru(p,q)$.
The $H\in i\frt_0$ corresponding to $\fru(p,q)\times\frso^*(2r)$ is
\[
H^{p, q}=(\unb{p}{1,\dots,1},\unb{r}{0,\dots,0},\unb{q}{-1,\dots,-1}).
\]
\end{prop}
\pf 
It is clear that the $H^{p,q}$ define admissible $\theta$-stable parabolic subalgebras with Levi subalgebras as in the statement of the proposition. We have to show that there are no other possibilities for a maximal Levi subalgebra.

Let $\frl$ be any Levi subalgebra of $\frg$ without compact factors defined by some $\frk$--dominant $H$. As before, we may assume that $H$ has integer coordinates. Since $H$ is $\frk$--dominant, $H$ is of the form
\[
H = (\unb{p_l}{l,\dots,l}, \unb{p_{l-1}}{l-1,\dots,l-1}, \dots, \unb{p_1}{1,\dots,1}, \unb{r}{0,\dots,0}, \unb{q_1}{-1,\dots,-1}, \dots, \unb{q_l}{-l,\dots,-l}),
\]
where $p_i, q_j \in \mathbb{Z}^{+}$. Since $\frl$ has no compact simple factors, 
\[
p_i=0\ \Rightarrow\ q_i=1;\qquad q_j=0\ \Rightarrow\ p_j=1.
\]
Obviously, 
\[
H^{'} = (\unb{p}{1,\dots,1}, \unb{r}{0,\dots,0}, \unb{q}{-1,\dots,-1}),
\] 
where $p = p_1 + p_2 + \cdots + p_l$, $q = q_1 + q_2 + \cdots + q_l$ defines a strictly larger $\frl$.

This proves (i) and (ii).

(iii) and (iv) are also proved, since we have seen that there are no further possibilities for $H$ corresponding to Levi subalgebra without compact factors besides the ones listed; in particular, each $\frl$ corresponds to only one $H$ and hence to only one $\frq$. 
\epf

Putting the results in Proposition \ref{leviD2} together with \cite[Proposition~6.1.]{HPV}, we obtain

\begin{cor}
\label{combparD2}
Let $\frg_0=\frso^*(2n)$ with $n \geq 4$. 
The modules $A_\frq$, or equivalently the admissible $\theta$-stable parabolic subalgebras $\frq=\frl\oplus\fru$ of $\frg$ such that $\frl_0$ has no compact factors, are in one-to-one correspondence with ordered sequences
\[
(p_1,q_1),(p_2,q_2),\dots,(p_l,q_l),
\]
with $l,p_i,q_j\in\bbZ^+$, such that
\begin{gather*}
p_1+p_2+\dots+p_l+q_1+q_2+\dots+q_l\leq n;\\
\text{and } \qquad p_i=0\ \Rightarrow\ q_i=1;\qquad q_j=0\ \Rightarrow\ p_j=1.
\end{gather*}
\qed
\end{cor}

Let us describe the faces of $W \rho$ corresponding to $A_{\frq}$ modules given by the proposition \ref{leviD2}. 

The shortest $w\in\wpg$ such that 
\[
w^{-1}H = w^{-1}(\unb{p_l}{l,\dots,l}, \unb{p_{l-1}}{l-1,\dots,l-1}, \dots \unb{p_1}{1,\dots,1}, \unb{r}{0,\dots,0}, \unb{q_1}{-1,\dots,-1}, \dots, \unb{q_l}{-l,\dots,-l})
\]
is $\frg$-dominant with respect to the standard $\dpg$ sends
$\rho=(n - 1,n - 2, \dots, 1, 0)$ to
\begin{align*}
\rho'=(& \unb{p_l}{n - 1, n - 2, \dots,n - p_l }, \unb{p_{l-1}}{n - (p_l + q_l) - 1,\dots,n - (p_l + q_l) - p_{l-1}}, \dots \\
& \unb{p_1}{r + p_1 + q_1 - 1, \dots, r + q_1}, \unb{r}{r - 1, r - 2, \dots, 0}, \unb{q_1}{- r, - (r +
1) \dots,-(r + q_1 - 1)}, \\
& \unb{q_2}{ - (r + p_1 + q_1), \dots,  - (r + p_1 + q_1 + q_2 - 1)}, \dots, \unb{q_l}{- n + p_l + q_l, \dots,  - n + p_l + 1}).
\end{align*}

The vertices of the corresponding face $\Phi = W_{\frl} \rho'$ are obtained from $\rho'$ by:
\begin{enumerate}
\item permuting coordinates in groups that correspond to repeated coordinates of $H$.
\item exchanging pairs $(a, -b)$ and $(b, -a)$ in groups of coordinates with the first component in the group that corresponds to repeated coordinates $t, \cdots, t$ of $H$, and the second component in the group that corresponds to repeated coordinates $-t, \cdots, -t$ of $H$, where $t \in \{1, \cdots, l\}$.
\item permuting and making an even number of sign changes to the group of coordinates on the places $p + 1, \cdots p + r$.
\end{enumerate}

Now $\Phi\cap C_\frk$, which corresponds to $H_D(\aq)$ for the $\frq$ corresponding to $\Phi$, consists of the vertices of $\Phi$ that are $\frk$-dominant, i.e., the vertices with coordinates in descending order. The descending order can be checked separately for each group of coordinates. From this one can see that each $\Phi$ is determined by $\Phi\cap C_\frk$, so the modules $\aq$ are determined by their Dirac cohomology. 

We summarize the results for the case $\frg_0=\frso^*(2n)$ in the following theorem.

\begin{thm}
\label{thmso*2n}
If $\frg_0=\frso^*(2n)$, then the modules $\aq$ correspond to faces $\Phi$ of $W\rho$ intersecting $C_\frk$, which can be described by the following procedure:

\begin{enumerate}
\item Pick a $\frk$-dominant conjugate $\rho'$ of $\rho$;
\item Pick an ordered sequence $(p_1,q_1),\dots,(p_l,q_l)$ as in Corollary \ref{combparD2}; let $p=p_1+\dots+p_l$ and $q=q_1+\dots+q_l$;
\item Consider the following groups of coordinates of $\rho'$:
\begin{multline*}
\rho'_1,\dots,\rho'_{p_1};\rho'_{n-q_l+1},\dots,\rho'_{n};\qquad \rho'_{p_1+1},\dots,\rho'_{p_1+p_2};\rho'_{n-q_l - q_{l-1} + 1},\dots,\rho'_{n-q_l};\qquad\dots \\
\dots\qquad \rho'_{p-p_l+1},\dots,\rho'_{p};\rho'_{n - q +1},\dots,\rho'_{n - q + q_1}.
\end{multline*}

Each of these groups of coordinates has two subgroups, separated by a semicolon. In the first of these subgroups, all coordinates are positive, and in the second subgroup the coordinates are negative. Now we first ignore the signs, permute the whole group in all possible ways, and then put back the signs to the coordinates in the new second subgroup.

The coordinates that are not part of the above groups of coordinates form the last group 
\[
\rho'_{p + 1},\dots,\rho'_{p+r}.
\]

To these coordinates we apply all permutations and even number of sign changes.
\end{enumerate}
\item If an $\aq$ module corresponds to a face $\Phi$, then
its Dirac cohomology $H_D(A_\frq)$ consists of the $\Kt$-types with infinitesimal characters in $\Phi\cap C_\frk$.
$H_D(\aq)$, or $\Phi\cap C_\frk$, uniquely determines the module $A_\frq$.  The set $\Phi\cap C_\frk$ consists of vertices that have decreasing coordinates, and this property can be checked for the whole vertex at once, or group by group. 
\qed
\end{thm}

\section{Classification of $A_{q}$ modules for type $G_2$}\label{sec:G2}
Let us now consider the case $\frg_0$ of type $G_2$. The unitary dual of $G_2$ is described in \cite{V3}. Here we use the notation from \cite{K}. The Weyl group is the dihedral group $D_6$. 
The standard positive roots are 
\[
\epsilon_1 - \epsilon_2, -\epsilon_1 + \epsilon_3, -\epsilon_2 + \epsilon_3, -2\epsilon_1 + \epsilon_2 + \epsilon_3, -\epsilon_1 - \epsilon_2 + 2 \epsilon_3, \epsilon_1 - 2\epsilon_2 + \epsilon_3.
\]
The maximal compact subalgebra of $\frg_0$ is $\frk_0 = \mathfrak{su}(2) \oplus \mathfrak{su}(2)$ (see \cite{K}) with simple roots 
\[
\epsilon_1 - \epsilon_2, -\epsilon_1 - \epsilon_2 + 2 \epsilon_3.
\]
We will denote the elements of $i\mathfrak{t}_0 \simeq \mathfrak{t}_{\mathbb{R}}^{*}$ by $(x_1, x_2, x_3)$, where $x_1 + x_2 + x_3 = 0$.
An element $H = (H_1, H_2, H_3) \in \mathfrak{t}_{\mathbb{R}}^{*}$ is $\frk$--dominant if and only if
\[
H_1 - H_2 \geq 0 \text{ and } -H_1 - H_2 + 2H_3 \geq 0, 
\]
that is, if and only if $$H_2 \leq H_1 \leq - H_2.$$
There are three positive root systems $\Delta^{+}$ containing the above system of compact positive roots and the corresponding half sums of positive roots are:
\[
\rho_1 = (-1, -2, 3), \quad \rho_2 = (1, -3, 2), \quad \rho_3 = (2, -3, 1).
\]
There are three discrete series representations of $G_2$ and each of them has $\Phi \cap C_{\frk}$ equal to one of $\rho_i$ with the corresponding $H = \rho_i$. In each of these cases the Levi subalgebra is $\mathfrak{l} = \mathfrak{t}$. The trivial module has $H = 0$ ($\mathfrak{l} = \mathfrak{g}$) and $\Phi \cap C_{\frk} = \{ \rho_1, \rho_2, \rho_3 \}$. The other Levi subalgebras of $\frg_0$ are
\begin{align*}
 \frl^1 & = \frt \oplus \frg_{-\epsilon_1 + \epsilon_3} \oplus \frg_{\epsilon_1 - \epsilon_3}, \quad & \frl^2 & = \frt \oplus \frg_{-2\epsilon_1 + \epsilon_2 + \epsilon_3} \oplus \frg_{2\epsilon_1 - \epsilon_2 - \epsilon_3}, 
\\ \frl^3 & = \frt \oplus \frg_{\epsilon_1 - \epsilon_2} \oplus \frg_{-\epsilon_1 + \epsilon_2}, \quad & \frl^4 & = \frt \oplus \frg_{-\epsilon_1 - \epsilon_2 + 2\epsilon_3} \oplus \frg_{\epsilon_1 + \epsilon_2 - 2\epsilon_3}.
\end{align*}
Among them $\frl^{1}_{0}$ and $\frl^{2}_{0}$ are without compact simple factors. The cooresponding $H$ are respectively $H_1 = (1, -2, 1)$ and $H_2 = (0, -1, 1)$. Since $H_2$ is $\frg$--dominant, the corresponding $\rho'$ is $\rho_1 = (-1, -2, 3)$, while the corresponding $\rho'$ for $H_1$ is $\rho_3 = (2, -3, 1)$. Now, it is easy to check that the corresponding faces are respectively
\begin{align*}
\Phi_1 & = W_{l^1} (2, -3, 1) = \{ (2, -3, 1), (1, -3, 2) \}, \\
\Phi_2 & = W_{l^2} (-1, -2, 3) = \{ (-1, -2, 3), (1, -3, 2) \}
\end{align*}
with the corresponding Dirac cohomologies
\[
\Phi_1 \cap C_{\frk} = \Phi_1 = \{ \rho_2, \rho_3\}, \quad \Phi_2 \cap C_{\frk} = \Phi_2 = \{ \rho_1, \rho_2 \}.
\]
We have proved the following theorem:
\begin{thm}
Let $\frg_0$ be the real form of $\frg_2$. Then there are three discrete series representations and each of them has $\Phi \cap C_{\frk}$ equal to one of $\rho_i$. The trivial module has the Dirac cohomology equal to $\Phi \cap C_{\frk} = \{ \rho_1, \rho_2, \rho_3 \}$. The other $A_{\frq}$ and their Dirac cohomologies are 
\begin{align*}
\Phi_1 & = \{ \rho_2, \rho_3\}, \quad  \Phi_1 \cap C_{\frk}  = \{ \rho_2, \rho_3\} \\
\Phi_2 & = \{ \rho_1, \rho_2 \}, \quad \Phi_2 \cap C_{\frk}  = \{ \rho_1, \rho_2 \}.
\end{align*}
All of the above $A_{\frq}$ modules are distinguished by their Dirac cohomology. 
\end{thm}

\end{document}